\documentclass[12pt]{article}
\usepackage{hyperref,url}
\usepackage{graphicx,amssymb}
\usepackage{amsmath,amsthm}
\usepackage{hyperref,url}
\usepackage{verbatim}
\usepackage{listings}
\newcommand{\nocopyright}{
No Copyright\thanks{
The author(s) hereby waive all copyright
and related or neighboring rights to this work,
and dedicate it to the public domain.
This applies worldwide.
}}
\title{Let's reinvent subtraction}
\author{Peter G. Doyle}
\date{Version 240 dated 25 March 2022\\
\nocopyright
}
\newtheorem{theorem}{Theorem}
\newtheorem{prop}[theorem]{Proposition}

\theoremstyle{definition}
\newtheorem*{definition}{Definition}

\def\proofstart{{\bf {\medskip}{\noindent}Proof: }}
\newcommand{\mathproofend}{\quad \spadesuit}
\newcommand{\proofend}{$\quad \spadesuit$}

\newcommand{\p}{\;\;\;}

\newcommand{\pd}{\p.}

\makeatletter
\newcommand{\xRrightarrow}[2][]{\ext@arrow 0359\Rrightarrowfill@{#1}{#2}}
\newcommand{\Rrightarrowfill@}{\arrowfill@\equiv\equiv\Rrightarrow}
\newcommand{\xLleftarrow}[2][]{\ext@arrow 3095\Lleftarrowfill@{#1}{#2}}
\newcommand{\Lleftarrowfill@}{\arrowfill@\Lleftarrow\equiv\equiv}
\newcommand{\xLleftRrightarrow}[2][]{\ext@arrow 3399\LleftRrightarrowfill@{#1}{#2}}
\newcommand{\LleftRrightarrowfill@}{\arrowfill@\Lleftarrow\equiv\Rrightarrow}
\makeatother

\newcommand{\cross}{\times}
\newcommand{\op}{\operatorname}

\newcommand{\without}[2]{#1 \setminus #2}

\newcommand{\kwstyle}{\tt}
\newcommand{\IF}{{\ \kwstyle{if}\ }}
\newcommand{\ELSE}{{\ \kwstyle{else}\ }}
\newcommand{\WHILE}{{\kwstyle{while}\ }}
\newcommand{\DO}{{\ \kwstyle{do}\ }}
\newcommand{\RETURN}{{\kwstyle{return}\ }}

\newcommand{\into}{\rightarrowtail}
\newcommand{\resp}{\ll}
\newcommand{\One}{\mathrm{One}}
\newcommand{\Dbar}{\bar{D}}
\newcommand{\fix}{\mathrm{fix}}

\newcommand{\xdiv}{\mathrm{xdiv}}
\begin{document}

\maketitle

\begin{abstract}
Subtraction is a powerful technique for creating new bijections from old.
Let's reinvent it!
While we're at it, let's reinvent division as well.
\end{abstract}

\section{Matchings}

Write
\[
f: A \equiv B
\]
and say `$f$ matches $A$ with $B$'
to mean that we know a suitable bijection
$f$ from $A$ to $B$,
together with its inverse
\[
f^{-1}:A \equiv B
.
\]
Write
\[
A \equiv B
\]
and say `$A$ matches $B$'
to mean that we know (or know we \emph{could} know)
some $f: A \equiv B$.
We have
\[
A \equiv A
;\p
A \equiv B \implies B \equiv A
;\p
A \equiv B \land B \equiv C \implies A \equiv C
.
\]
(We refrain from saying that $\equiv$ is an equivalence relation
since it
is inherently time-dependent.)

We can add and multiply matchings:
\[
A \equiv B \land C \equiv D 
\implies
A + C \equiv B + D
,
\]
where $+$ denotes disjoint union, and
\[
A \equiv B \land C \equiv D 
\implies
A \cross C \equiv B \cross D
.
\]

\section{Respectful subtraction}

Before addressing subtraction in general, let's begin with
\emph{respectful subtraction}, an important special case.
It's so simple that it hardly deserves to be called subtraction.

\begin{definition}
For $f: A+C \equiv B+D$, $g:C \equiv D$, write
$g \resp f$ and say `$g$ respects $f$' if
\[
\forall x \in C
\;
(f(x) \in D \implies g(x)=f(x))
.
\]
\end{definition}

\begin{prop}[Respectful subtraction]
If
\[
f: A+C \equiv B+D
;\p
g:C \equiv D
;\p
g \resp f
\]
then
\[
\without{f}{g}: A \equiv B
,
\]
where
\[
\without{f}{g} (x) =
f(x) \IF f(x) \in B \ELSE f(g^{-1}(f(x)))
.
\]
Moreover,
\[
\without{f}{g} \resp f
\]
and
\[
\without{f}({\without{f}{g}}) = g
.
\]
\end{prop}

\proofstart
Without loss of generality,
assume $f$ is the identity on $A+B=C+D$.
$g$ fixes $C \cap D$ and matches $C \setminus D$ to $D \setminus C$.
$\without{f}{g}$ fixes $A \cap B$ and, taking its cue from $g^{-1}$,
matches $A \setminus B=D \setminus C$ to $B \setminus A=C \setminus D$.
\proofend

\section{Subtraction}

\begin{prop}[Subtraction] \label{subtraction}
If
\[
f: A+C \equiv B+D
;\p
g:C \equiv D
\]
with $D$ finite,
then
\[
\without{f}{g}: A \equiv B
,
\]
where
\[
\without{f}{g}(x) =
(
y:=f(x);\;
\WHILE y \in D
\DO y:=f(g^{-1}(y));\;
\RETURN y
)
.
\]
\end{prop}

\proofstart
This goes way back---see \cite{doyle:category}.
\proofend
\\

When $g \resp f$ we're back to respectful subtraction:
\begin{prop}
If 
$g \resp f$ then
\[
\without{f}{g} (x) =
f(x) \IF f(x) \in B \ELSE f(g^{-1}(f(x)))
.
\mathproofend
\]
\end{prop}

The fact that $\without{f}{g} \resp f$ is general:
\begin{prop} \label{prop2}
$\without{f}{g} \resp f$.
\proofend
\end{prop}

Idempotence of subtraction characterizes respectfulness:

\begin{prop} \label{prop3}
\[
\without{f}{(\without{f}{g})} = g
\;\iff\;
g \resp f
.
\mathproofend
\]
\end{prop}

This gives us a closure operation:

\begin{prop} \label{prop4}
\[
\without{f}{(\without{f}{(\without{f}{g})})} = \without{f}{g}
,
\]
so
\[
\without{f}{(\without{f}{(\without{f}{(\without{f}{g})})})}
= \without{f}{(\without{f}{g})}
.
\mathproofend
\]
\end{prop}

\section{Inclusion-exclusion}

Subtraction generalizes to inclusion-exclusion.
Let $P$ be a poset with
$\{q:q\leq p \}$ finite for all $p$.
Given a family of finite sets $A_p, p \in P$, let
\[
A_{\leq p} = \sum_{q \leq p} A_q
,
\]
etc.

\begin{prop}[Inclusion-exclusion] \label{incexc}
\[
\forall p
\;
A_{\leq p} \equiv B_{\leq p}
\;\;\implies\;\;
\forall p
\;
A_p \equiv B_p
.
\]
\end{prop}

\proofstart
By induction:
Assuming
\[
\forall q <p
\;
A_q \equiv B_q
\]
(true if $p$ is minimal)
we have
\[
A_{<p} \equiv B_{<p}
.
\]
Subtract from
\[
A_{\leq p} \equiv B_{\leq p}
\]
to get
\[
A_p \equiv B_p
.
\mathproofend
\]

To be more explicit,
define the projection map
\[
\pi_A: \sum_p A_p \to P
,
\;
\pi(x) = p \iff x \in A_p
.
\]

\begin{prop}
If
\[
g_p: A_{\leq p} \equiv B_{\leq p}
\]
then
\[
f_p: A_p \equiv B_p
\]
where
\[
f_p(x)=F(p,x)
,
\]
\[
F(p,x) =
(
y:=g_p(x);\;
q:=\pi_B(y);\;
\RETURN
y \IF q=p
\ELSE
F(p,\bar{F}(q,y))
)
;\;
\]
\[
\bar{F}(p,x) =
(
y:=g_p^{-1}(x);\;
q:=\pi_A(y);\;
\RETURN
y \IF q=p
\ELSE
\bar{F}(p,F(q,y))
)
.
\]
\end{prop}

\proofstart
This is what you get if you trace it through.
\proofend

\section{Extreme division}

Subtraction will get you a long way toward automating the process
of generating bijections.
But sometimes you will want to divide, and that's when
things can get scary.
(Like the old Marchant mechanical calculators, which would make
a satisfying `chunk' when you hit $+$ or $-$,
but would make a terrifying racket, with the carriage scurrying to and fro,
when you hit the ${\overset{\mbox{\tiny Auto}}{\div}}$ key.)

The way to keep things under control is to make sure you are multiplying
polynomials (in one or many variables),
with the polynomial you are dividing
by having a unique extreme monomial $\omega$
for some linear function on the space of
degrees (in other words a singleton monomial on the boundary of its
Newton polytope).
In this case you can use `extreme division', whereby you recursively
subtract the mapping based on multiplication by $\omega$.

Suppose
\[
F:A \cross C \equiv B \cross C
\]
and
\[
G:B \cross C \equiv A \cross C
\pd
\]
(We may choose $G=F^{-1}$, but we don't require this.)
For $\omega \in C$,
define
\[
\op{xdiv}((F,G),\omega)
=
(f,g)
\]
where $f,g$ are the
partial functions $f$ on $A$ and $g$ on $B$
defined via the mutual recursion equations
\[
f(x)
=
((y,z):=F((x,\omega))
;\,
\WHILE
z \neq \omega
\DO
((y,z):=F((g(y),z))
;\,
\RETURN
y)
\]
\[
g(x)
=
((y,z):=G((x,\omega))
;\,
\WHILE
z \neq \omega
\DO
((y,z):=G((f(y),z))
;\,
\RETURN
y)
\]

If both $f$ and $g$ are total, we say that
the pair $(F,G)$ is
\emph{X-divisible for $\omega$},
and that $\omega$ is an \emph{extreme point} for the pair $(F,G)$.
This terminology springs from the following proposition.

\begin{prop}[Extreme division] \label{xdiv}
If $A,B,C$ are multinomials, and $F,G$ match terms of
$A \cdot C$ and $B \cdot C$,
then $(F,G)$ is X-divisible for any extreme monomial $\omega$ of $C$.
\end{prop}

\proofstart
By induction.
\proofend

Extreme division is what Conway and Doyle
\cite{conwaydoyle:three}
had in mind when they wrote,
`There is more to division than repeated subtraction.'
This must be a great truth, because its negation
would appear to be at least as true.

\section{Mode d'emploi}

Contrary to what we may appear to be claiming in Proposition \ref{xdiv},
the bijection yielded by the $\xdiv$ algorithm may fail to be a `matching',
because it may take too long to compute,
or otherwise fail to qualify as `suitable',
the admittedly slippery condition that we slipped in as
part of the definition of a matching.
The same goes for bijections obtained by inclusion-exclusion.
This is why we have taken care to announce \ref{incexc} and \ref{xdiv}
as `Propositions', yielding bijections
proposed for consideration as `matchings'.
This is contrary to mathematical custom,
and wrong-headed,
but useful nevertheless.

Now it will often happen that a slow quotient bijection can be speeded up
immensely by `memoizing', meaning that values
of $f$ and $g$ are automatically saved so that they don't get 
computed over and over.
This in itself may make the bijection `suitable'.

Better still is to be able to see `what the bijection is doing',
so that it can be defined, and proven to be a suitable bijection,
without reference to its origin as a quotient.

Here's a case in point.
Everyone knows that
\[
\binom{n}{k}
=
\binom{n}{n-k}
,
\]
and if you ask why,
they will either tell you to match a $k$-subset to
the complementary $(n-k)$ - subset,
or compute
\[
\binom{n}{k}
=
\frac{n!}{k!(n-k)!}
=
\frac{n!}{(n-k)!(n-(n-k))!}
=
\binom{n}{n-k}
.
\]

Taking our motto to be `follow the algebra',
we recast the computation as 
\[
\binom{n}{k} k! (n-k)!
=
n!
=
\binom{n}{n-k} (n-k)! k!
=
\binom{n}{n-k} k! (n-k)!
,
\]
where every step is backed by a matching.
Now divide.
The xdiv algorithm yields a very inefficient computation of a very simple
matching
(see the code in the appendix):
\clearpage
\begin{verbatim}
([[0, 1], [2, 3, 4]], [[0, 1, 2], [3, 4]])
([[0, 2], [1, 3, 4]], [[0, 1, 3], [2, 4]])
([[0, 3], [1, 2, 4]], [[0, 2, 3], [1, 4]])
([[0, 4], [1, 2, 3]], [[1, 2, 3], [0, 4]])
([[1, 2], [0, 3, 4]], [[0, 1, 4], [2, 3]])
([[1, 3], [0, 2, 4]], [[0, 2, 4], [1, 3]])
([[1, 4], [0, 2, 3]], [[1, 2, 4], [0, 3]])
([[2, 3], [0, 1, 4]], [[0, 3, 4], [1, 2]])
([[2, 4], [0, 1, 3]], [[1, 3, 4], [0, 2]])
([[3, 4], [0, 1, 2]], [[2, 3, 4], [0, 1]])
\end{verbatim}
This follow-the-algebra matching
differs from only slightly from taking the complementary set.
It's arguably better.
Do you agree?

\clearpage

\section*{Appendix}

\begin{verbatim}
"""
gauss 000 - use extreme division to match (n choose k) to (n choose n-k)
"""

def xdiv(FG,omega):
  F,G=FG
  def f(x):
    y,z=F((x,omega))
    while z!=omega:
      y,z=F((g(y),z))
    return y
  def g(x):
    y,z=G((x,omega))
    while z!=omega:
      y,z=G((f(y),z))
    return y
  return [f,g]

import itertools

def makeintolist(a): return [x for x in a]
def sublist(a,c): return [a[x] for x in c]
def num(n): return [x for x in range(n)]

def combinations(a,k):
  return [makeintolist(a) for a in itertools.combinations(a,k)]

def complement(s,a): return [x for x in s if x not in a]
def separate(s,a): return [a,complement(s,a)]

"""
choose(n,k) lists the ways for separating n into pieces of size k and n-k 
"""

def choose(n,k):
  s=range(n)
  return [separate(s,a) for a in combinations(s,k)]


"""
compress(l) replaces elements of l by their relative ranks
"""

def compress(l):
  m=sorted(l)
  return [m.index(x) for x in l]

"""
binom(sigma,k) maps the permutation sigma to choose(n,k) * k! * (n-k)!
"""

def binom(sigma,k):
  a=sigma[:k]
  b=sigma[k:]
  return [[sorted(a),sorted(b)],[compress(a),compress(b)]]

"""
monib is the inverse of binom
"""

def monib(abcd):
  ((a,b),(c,d))=abcd
  return sublist(a,c)+sublist(b,d)

"""
flipkl maps  choose(n,k) * k! * (n-k)! to  choose(n,k) * (n-k)! * k!

We need this because we decided to cast division in terms of matchings
between A*C and B*C rather than A*C and B*D
"""

def flipkl(abcd):
  ((a,b),(c,d))=abcd
  return [[a,b],[d,c]]

"""
Our F and G would be the same, except for the flipping.
"""

def F(abcd):
  ((a,b),(c,d))=abcd
  k=len(a)
  l=len(b)
  return flipkl(binom(monib(abcd),l))

def G(abcd):
  ((a,b),(c,d))=abcd
  k=len(a)
  l=len(b)
  return binom(monib(flipkl(abcd)),l)

def omega(n,k): return [num(k),num(n-k)]

def match(n,k): return xdiv([F,G],omega(n,k))

"""
TESTING
"""

def column(l):
  print(*l,sep='\n')

n=5
k=2

A=choose(n,k)
(f,g)=match(n,k)
B=[f(ab) for ab in A]
column(zip(A,B))
C=[g(ab) for ab in B]
print(A==C)
\end{verbatim}

\bibliography{lrs}
\bibliographystyle{hplain}
\end{document}